\documentclass[a4paper,12pt]{article}
\usepackage{amsmath,amsthm,amssymb,latexsym}
\title{{\bf The capacity associated to signed Riesz kernels, and Wolff potentials.}}
\author{\Large{\Large Joan Mateu, Laura Prat and Joan Verdera.}}
\begin{document}
\date{}
\maketitle
%\date
\newtheorem{teo}{Theorem}[section]
\newtheorem{co}[teo]{Corollary}
\newtheorem{lemma}[teo]{Lemma}
\newtheorem{defi}[teo]{Definition}
\newtheorem{note}[teo]{Note}
\newcommand{\Ha}{{\cal H}^{\alpha}}
\newcommand{\Rn}{{\mathbb R}^n}
\newcommand{\ep}{\varepsilon}
\newcommand{\N}{\mathbb{N}}
\newcommand{\Z}{\mathbb{Z}}
\newcommand{\C}{\mathbb{C}}

\begin{abstract}
We show that, for $0<\alpha<1$, the capacity associated to the signed vector valued
Riesz kernel $\frac x{|x|^{1+\alpha}}$ in $\Rn$ is comparable to the Riesz capacity
$C_{\frac 2 3(n-\alpha),\frac 3 2}$ of non-linear potential theory.
\end{abstract}

\section {Introduction.}

In this paper we study the capacity $\gamma_\alpha$ associated to the signed
 vector valued Riesz kernels
$k_\alpha(x)=\frac{x}{|x|^{1+\alpha}}$, $0<\alpha<n$, in $\Rn$. If $K\subset\Rn$ is
compact one sets

$$\gamma_{\alpha}(K)=\sup|<T,1>|,$$

\noindent where the supremum is taken over all distributions $T$ supported on $K$ such
that $\;T*\frac{x_i}{|x|^{1+\alpha}}$ is a function in $\;L^{\infty}(\mathbb{R}^n)$
and $\|T*\frac{x_i}{|x|^{1+\alpha}}\|_{\infty}\leq 1$, for $\;1\leq i\leq n$. For
$n=2$ and $\alpha=1$ this is basically analytic capacity (see \cite{semiadditivity}),
and for $\alpha=n-1$ and any $n\geq 2$, $\gamma_{n-1}$ is Lipschitz harmonic capacity
(see \cite{paramonov}, \cite{mp} and \cite{vecm}).

In \cite{laura} one discovered the fact that if $0<\alpha<1$, then a compact set of
finite $\alpha$-dimensional Hausdorff measure has zero $\gamma_\alpha$ capacity. This
is in strong contrast with the situation for integer $\alpha$, in which
$\alpha$-dimensional smooth hypersurfaces have positive $\gamma_\alpha$ capacity. The
case of non-integer $\alpha>1$ is not completely understood, although it was shown in
\cite{laura} that for Ahlfors-David regular sets the result mentioned above for
$0<\alpha<1$ still holds.

In this paper we establish the equivalence between $\gamma_\alpha$, $0<\alpha<1$, and
one of the well-known Riesz capacities of non-linear potential theory (see
\cite{adamshedbergllibre}, Chapter 1, p. 38). The Riesz capacity $C_{s,p}$ of a
compact set $K\subset\Rn$, $1<p<\infty$, $0<sp\leq n$, is defined by

$$C_{s,p}(K)=\inf\{\|\varphi\|_p^p:\;\varphi*\frac 1{|x|^{n-s}}\geq 1\;\mbox{ on }K\},$$

\noindent where the infimum is taken over all compactly supported infinitely
differentiable functions on $\Rn$. The capacity $C_{s,p}$ plays a central role in
understanding the nature of Sobolev spaces (see \cite{adamshedbergllibre}).

Our main result is the following surprising inequality.\newline

{\bf Theorem.} {\em For each compact set $K\subset\Rn$ and for $0<\alpha<1$ we have

$$C^{-1}C_{\frac 2 3(n-\alpha), \frac 3
2}(K)\leq\gamma_\alpha(K)\leq C\;C_{\frac 2 3(n-\alpha),\frac 3 2}(K),$$ where $C$ is
a positive constant depending only on $\alpha$ and $n$.}

\vspace{.5cm} Since it is well-known that $C_{\frac 2 3(n-\alpha), \frac 3 2}$
vanishes on sets of finite $\alpha$-dimensional Hausdorff measure (see
\cite{adamshedbergllibre}, Theorem 5.1.9, p.134), the same applies to $\gamma_\alpha$.
Thus we recover one of the main results of \cite{laura}. On the other hand, $C_{s,p}$
is a subadditive set function (almost by definition, see \cite{adamshedbergllibre},
p.26), and consequently, $\gamma_\alpha$ is semiadditive for $0<\alpha<1$, that is,
given compact sets $K_1$ and $K_2$,

\begin{equation}
\label{semia} \gamma_\alpha(K_1\cup K_2)\leq
C\left\{\gamma_\alpha(K_1)+\gamma_\alpha(K_2)\right\},
\end{equation}

\noindent for some constant $C$ depending only on $\alpha$ and $n$. In fact
$\gamma_\alpha$ is countably semiadditive. For $\alpha=1$ and $n=2$ inequality
(\ref{semia}) is still true and is a remarkable result obtained in
\cite{semiadditivity}. For $\alpha=n-1$ and any $n$ (\ref{semia}) has been shown very
recently in \cite{lipschvolberg}.\newline

Another interesting consequence of the Theorem is that $\gamma_\alpha$ is a
bilipschitz invariant. This means that if $\phi:\Rn\to\Rn$ is a bilipschitz
homeomorphism of $\Rn$, namely,

$$L^{-1}|x-y|\leq|\phi(x)-\phi(y)|\leq L|x-y|,$$

\noindent for $x,\;y\in\Rn$ and for some constant $L>0$, then for compact sets $K$ one
has

$$C^{-1}\gamma_\alpha(K)\leq\gamma_\alpha(\phi(K))\leq C\gamma_\alpha(K),$$

\noindent where $C$ depends only on $L$, $\alpha$ and $n$.

The bilipschitz invariance of the analytic capacity $\gamma$ has been recently proved
by X. Tolsa (see \cite{bilipschitz}). The result for a big class of Cantor sets was
proved before by Garnett and Verdera (see \cite{garnettverdera}).

Volberg has pointed out to the authors that a particular instance of the Theorem gives
the following curious result about Cauchy integrals. Take $n=2$ and
$\alpha=\frac{1}{2}.$ Then, given a compact set $K\subset \mathbb C,$ there exists a
distribution $T\neq 0$ supported on $K$ such that $T*\frac{z}{|z|^{3/2}}\in L^\infty
(\mathbb C)$ if and only if there exists a probability measure $\mu$ supported on $K$
such that $\mu*\frac{1}{z}\in L3(\mathbb C).$ This follows from the dual definition of
$C_{1,\frac{3}{2}}$ (see \cite{adamshedbergllibre}, Theorem 2.2.7.).

Our proof of the Theorem rests on two steps. The first one is the analogue for
$0<\alpha<1$ of the main result in \cite{semiadditivity}, namely, the equivalence
between $\gamma_\alpha$ and $\gamma_{\alpha,+}$. For a compact set $K\subset\Rn$, the
positive $\gamma_\alpha$ capacity is defined by

$$\gamma_{\alpha,+}(K)=\sup\mu(K),$$

\noindent where the supremum is taken over those positive Radon measures $\mu$
supported on $K$ such that $\frac{x_i}{|x|^{1+\alpha}}*\mu$ is in $L^\infty(\Rn)$ and
$\left\|\frac{x_i}{|x|^{1+\alpha}}*\mu\right\|_\infty\leq 1$, for $1\leq i\leq n$.
Clearly $\gamma_{\alpha,+}(K)\leq \gamma_{\alpha}(K)$ for any $K.$
\begin{teo}
\label{semiadmeu} For each compact set $K\subset\Rn$ and $0<\alpha<1$, we have

$$\gamma_{\alpha,+}(K)\leq\gamma_\alpha(K)\leq
C\gamma_{\alpha,+}(K),$$

\noindent where $C$ is some positive constant depending only on $\alpha$ and $n$.
\end{teo}

We claim that Theorem \ref{semiadmeu} can be proved by adapting the scheme of the
proof of Theorem 1.1 in \cite{semiadditivity} and the adjustments introduced in
\cite{semiadditivity2} to prove Theorem 7.1 there. This is explained in some detail in
section \ref{section2}. When analyzing the argument used in \cite{semiadditivity} one
realizes that it is based on two main technical ingredients. The first is the
non-negativity of the quantity obtained when symmetrizing the kernel, which was proved
in \cite{laura} for the Riesz kernel $k_\alpha$ with $0<\alpha<1$. The second is the
fact that the Cauchy kernel ( that is, $k_1$ in dimension $n=2$) localizes in the
uniform norm. By this we mean that if $T$ is a compactly supported distribution such
that $T*\frac{1}{z}$ is a bounded function then $(\varphi T)*\frac{1}{z}$ is also
bounded for each compactly supported ${\cal C}^1$ function $\varphi$ and we have the
corresponding estimate. This is an old result, which is simple to prove because
$\frac{1}{z}$ is related to the differential operator $\overset-
\partial$ (\cite{garnett}, Chapter V). The same localization
result can be proved easily for any $n$ and $\alpha=n-1,$ because $k_{n-1}$ is related
to the Laplacian (\cite{paramonov} and \cite{vecm}). For other parameters $\alpha$
between $0$ and $n$ is not clear at all that there is a differential operator in the
background and consequently the corresponding localization result becomes far from
being obvious. In fact, the proof of the localization Theorem for $k_{\alpha}$ for any
$\alpha, \, 0<\alpha<n,$ is the main technical obstacle we have to surmount in this
paper. When localization is available there is no obstruction in adapting Lemma 7.2
(part $(h)$) in \cite{semiadditivity2}. Once Theorem \ref{semiadmeu} is at our
disposal we need to relate $\gamma_{\alpha,+}$ to $C_{\frac 2 3(n-\alpha),\frac 3 2}$
and this is the second step in the proof of the Theorem. \newline

The plan of the paper is the following. Section 2 contains some preliminary
definitions and results that will be used throughout the article. In section 3 we
prove the localization theorem for the signed Riesz potentials. In section 4 we
complete the proof of the main Theorem showing that $\gamma_{\alpha,+}$ is comparable
to $C_{\frac 2 3(n-\alpha), \frac 3 2}.$
\newline

Constants independent of the relevant parameters are denoted by $C$ and may be
different at each occurrence. The notation $A\approx B$ means, as it is usual, that
for some constant $C$ one has $C^{-1}B\leq A\leq C B$.

\section{Preliminaries.}
\subsection {Simmetrization of Riesz kernels.}
The symmetrization process for the Cauchy kernel introduced in \cite{melnikovdiscret}
has been successfully applied in these last years to many problems of analytic
capacity and $L2$ boundedness of the Cauchy integral operator (see \cite{mv},
\cite{mmv} for example; the survey \cite{davidpubmat} and the book \cite{pajotllibre}
contain many other interesting references). Given 3 distinct points in the plane,
$z_1, z_2$ and $z_3$, one finds out, by an elementary computation that

\begin{equation}
c(z_1,z_2,z_3)2=\sum_{\sigma}\frac{1}{(z_{\sigma(1)}-z_{\sigma(3)})
\overline{(z_{\sigma(2)}-z_{\sigma(3)}})} \label{id}
\end{equation}

\noindent where the sum is taken over the six permutations of the set $\{1,2,3\}$ and
$c(z_1,z_2,z_3)$ is Menger curvature, that is, the inverse of the radius of the circle
through $z_1, z_2$ and $z_3$. In particular (\ref{id}) shows that the sum on the right
hand side is a non-negative quantity.

It can be shown that for $0<\alpha<1$ the symmetrization of the Riesz kernel
$k_{\alpha}(x)= x/|x|^{1+\alpha}$, gives also a positive quantity. On the other hand,
for $1<\alpha<n$, the phenomenon of change of signs appears when symmetrizing the
kernel $k_\alpha$, as one can easily check.

For $0<\alpha<n$ the quantity

\begin{equation} \label{per}
\sum_{\sigma}\frac{x_{\sigma(2)}-x_{\sigma(1)}}{|x_{\sigma(2)}-x_{\sigma(1)}|^{1+\alpha}}
\frac{x_{\sigma(3)}-x_{\sigma(1)}}{|x_{\sigma(3)}-x_{\sigma(1)}|^{1+\alpha}},
\end{equation}

\noindent where the sum is taken over the six permutations of the set $\{1,2,3\}$, is
the obvious analogue of the right hand side of (\ref{id}) for the Riesz kernel
$k_\alpha$. Notice that (\ref{per}) is exactly

$$ 2\;p_{\alpha}(x_1,x_2,x_3),$$

\noindent where $p_{\alpha}(x_1,x_2,x_3)$ is defined as the sum in (\ref{per}) taken
only on the three permutations $(1,2,3),\;(2,3,1), \; (3,1,2)$.\newline

In the following lemma we state the explicit description that was found in
\cite{laura} for the symmetrization of the Riesz kernel $k_\alpha$, for $0<\alpha <1$.

\begin{lemma}\label{acotaci}
Let $0<\alpha<1$, and $x_1,\;x_2,\;x_3$ three distinct points in $\mathbb{R}^n$.  Then
we have

$$\frac{2-2^\alpha}{L(x_1,x_2,x_3)^{2\alpha}}\leq p_{\alpha}(x_1,x_2,x_3)\leq
\frac{2^{1+\alpha}}{L(x_1,x_2,x_3)^{2\alpha}},$$

\noindent where $L(x_1,x_2,x_3)$ is the largest side of the triangle determined by
$x_1,\;x_2$ \; and \;$x_3$. In particular $p_{\alpha}(x_1,x_2,x_3)$ is a positive
quantity.
\end{lemma}

The relationship between the quantity $p_\alpha(x,y,z)$ and the $L2$ estimates of the
operator with kernel $k_\alpha$ is as follows. Take a positive finite Radon measure
$\mu$ in $\Rn$ which satisfies the growth condition $\mu(B(x,r))\leq r^\alpha$,
$x\in\Rn$, $r>0$. Given $\ep>0$, set

$$R_{\alpha,\ep}(\mu)(x)=\int_{|y-x|>\ep}k_\alpha(y-x)d\mu(y).$$

Then (see in \cite{mv} or \cite{pajotllibre} the argument for $\alpha=1$)

$$\left|\int\left|R_{\alpha,\ep}(\mu)(x)\right|^2d\mu(x)-\frac 1 3 p_{\alpha,\ep}(\mu)\right|\leq C\|\mu\|,$$

\noindent where $C$ is a constant depending only on $\alpha$ and $n$, and

$$p_{\alpha,\ep}(\mu)=\underset{S_\ep}{\iiint}p_\alpha(x,y,z)d\mu(x)d\mu(y)d\mu(z),$$

\noindent with

$$S_\ep=\{(x,y,z):\;|x-y|>\ep,\;|x-z|>\ep \mbox{ and }|y-z|>\ep\}.$$

Thus

\begin{equation}
\label{inequa} p_\alpha(\mu)\leq
3\sup_{\ep>0}\int\left|R_{\alpha,\ep}(\mu)(x)\right|^2d\mu(x)+C\|\mu\|,
\end{equation}

\noindent where

$$p_\alpha(\mu)=\int_{\Rn}\;\int_{\Rn}\;\int_{\Rn}p_\alpha(x,y,z)d\mu(x)d\mu(y)d\mu(z).\vspace{.5cm}$$

\subsection{The scheme of the proof of Theorem \ref{semiadmeu}.}
\label{section2} In this section we give an outline of the arguments involved in the
proof of Theorem \ref{semiadmeu}. The proof uses an induction argument on scales,
analogous to the one in \cite{mtv} and \cite{semiadditivity}. The main idea is to
show, by induction, that

$$\gamma_{\alpha,+}(K\cap Q)\approx\gamma_\alpha(K\cap Q)$$

\noindent for squares $Q$ of any size.\newline

The starting point in the proof of Theorem 1.1 in \cite{semiadditivity} is the
construction of a special family of cubes  $\{Q_j\}_{j=1}^N$ that cover $K$ and
satisfy

$$\gamma_{\alpha,+}(\cup_{j=1}^N Q_j)\leq C\gamma_{\alpha,+}(K)$$

\noindent and

$$\sum_{j=1}^N\gamma_{\alpha,+}(3Q_j\cap K)\leq
C\gamma_{\alpha,+}(K).$$

The construction of these cubes works without difficulty in the same way as in
\cite{semiadditivity} for $0<\alpha<1$, because we have non-negativity of the quantity
obtained when symmetrizing the Riesz kernel (see Lemma \ref{acotaci} above).\newline

>From the definition of the capacity $\gamma_\alpha$, it follows that
there exists a distribution $T_0$ supported on $K$ such that

\begin{enumerate}
\item$\displaystyle{\gamma_\alpha(K)\geq\frac 1 2\left|\left<T_0,1\right>\right|},$
\item $\displaystyle{\|T_0*\frac{x_i}{|x|^{1+\alpha}}\|_\infty\leq 1,\;\;1\leq i\leq n.}$
\end{enumerate}

Consider now a family of infinitely differentiable functions $\{\varphi_j\}_{j=1}^N$
such that each $\varphi_j$ is compactly supported on $2Q_j$, $0\leq\varphi_j\leq 1$,
$\displaystyle{\|\partial^s\varphi_j\|_\infty\leq\frac{C}{\ell(Q_j)^{|s|}}}\;$, $0\leq
|s|\leq n$ , and $\;\sum_{j=1}^N \varphi_j=1$   on $\cup_{j=1}^N Q_j$. At this point
we need an inequality of the type

$$\|\varphi_j T_0*\frac{x_i}{|x|^{1+\alpha}}\|_\infty\leq C$$

\noindent for $1\leq i\leq n$, $1\leq j\leq N$ and $0<\alpha<n$, with $C=C(\alpha,n).$
This will be proved in section \ref{section3}. Then, by definition of $\gamma_\alpha$
, we will obtain that

\begin{equation}
\label{capacitatlocal} \left|\left<\varphi_jT_0,1\right>\right|\leq
C\gamma_{\alpha}(2Q_j\cap K).
\end{equation}

\noindent for $1\leq j\leq N$.

Inequality (\ref{capacitatlocal}) is used later on in the proof in order to construct
a bounded function $b$ to which a suitable variant of the $T(b)$ theorem will be
applied. There is still one more difficulty in applying the Nazarov, Treil and Volberg
$T(b)$-type theorem one needs, namely, finding a substitute for what they call the
suppressed operators. It was already explained in \cite{laura} that there are at least
two versions of such operators for the Riesz kernels that work appropriately.

\section{Localization of Riesz potentials.}
\label{section3} One of the ingredients of the proof of Theorem 1.1 in
\cite{semiadditivity} is the localization of the Cauchy potential. The localization
method for the Cauchy potential, $T*1/z$, developed by A.G. Vitushkin for rational
approximation was adapted in \cite{paramonov} to localize the potential $T*x/|x|^n$
and used in problems of ${\cal C}^1$-harmonic approximation.\newline In this section
we will be concerned with the localization of the vector valued $\alpha-$Riesz
potentials $T*x/|x|^{1+\alpha}$, $0<\alpha<n$. \newline

Let $x=(x_1,...,x_n)\in\Rn$ and $|x|=\left(\sum_{i=1}^n x_i2\right)^{1/2}$. For
$s=(s_1,...,s_n)$, $\;0\leq s_i\in\mathbb{Z}\;$, we set
$x^s=x_1^{s_1}\cdots\;x_n^{s_n}$, $s!=s_1!\cdots\;s_n!$, $|s|=s_1+s_2+\cdots+s_n$,
$\partial^s=\partial^{s_1}/\partial x_1^{s_1}\cdots\;\partial^{s_n} /\partial
x_n^{s_n}$, $\Delta=\sum_{i=1}^n\partial2/\partial x_i2$ and
$\partial_j=\partial/\partial x_j$, $1\leq j\leq n$. In what follows, given a cube
$Q\subset\Rn$, $\varphi_Q$ will denote an infinitely differentiable function supported
on $2Q$ and such that $\|\partial^s\varphi_Q\|_\infty\leq C_s\ell(Q)^{-|s|}$, $0\leq
|s|\leq n$.\newline

\noindent We prove now the following general localization lemma.

\begin{lemma} \label{localitzacio} Let $0<\alpha<n$ and let $T$
be a compactly supported distribution such that $T*\frac{x_i}{|x|^{1+\alpha}}$ is a
bounded measurable function for $1\leq i\leq n$.
 Then there exists some constant $C=C(n,\alpha)>0$ such that

 $$\sup_{1\leq i\leq n} \|\varphi_Q T*\frac{x_i}{|x|^{1+\alpha}}\|_\infty
 \leq C\sup_{1\leq i\leq n}\|T*\frac{x_i}{|x|^{1+\alpha}}\|_\infty.$$
\end{lemma}

{\em Proof.} Our argument uses a reproduction formula for test functions involving the
kernel $\displaystyle{k^i(y)=\frac{y_i}{|y|^{1+\alpha}}},$ which was first introduced
in \cite{laura} (see Lemma 11). There are many variants of this formula depending, for
instance, on whether the dimension $n$ and the integer part of $\alpha$ are even or
odd. We will consider in full detail only the case of odd dimension of the form
$n=2k+1.$ We will also assume that $\alpha$ is non-integer and that its integer part
is even, of the form $[\alpha]=2d.$ At the end of the proof we shall briefly indicate
how to treat the remaining cases, including the case of integer $\alpha.$

Fix $x\in\mathbb R^n$ and set
 $$k^i_x(y)=\frac{x_i-y_i}{|x-y|^{1+\alpha}}.$$

We distinguish two cases:

\begin{enumerate}

\item[Case 1:] $x\in (3Q)^c$. Set $g(y)=\varphi_Q(y)k^i_x(y)$.
Lemma 11 in \cite{laura} tells us that

\begin{equation}
\label{representationformula}
 g(x)=c_{n,\alpha}\sum_{j=1}^n\left(\Delta^k\partial_jg*\frac{1}{|y|^{n-\alpha}}*k^j\right)(x),
\end{equation}

\noindent for some constant $c_{n,\alpha}$ depending only on $n$ and $\alpha$. We
emphasize that (\ref{representationformula}) works because $n$ is odd. Thus

$$\begin{array}{l}
\displaystyle{\left(\varphi_QT*k^i\right)(x)=<T,g>=c_{n,\alpha}\sum_{j=1}^n<
T*k^j,\Delta^k\partial_j g*\frac {1}{|y|^{n-\alpha}}>,}
\end{array}$$

and so

\begin{equation}
\begin{array}{l}
\label{primerpas} \displaystyle{\left(\varphi_Q
T*k^i\right)(x)=\sum_{j=1}^nc_{n,\alpha}\int_{(3Q)^c}(T*k^j)(z)\left(\Delta^k\partial_j
g*\frac {1}{|y|^{n-\alpha}}\right)(z)dz}\\\\
\displaystyle{+\sum_{j=1}^nc_{n,\alpha}\int_{3Q}(T*k^j)(z)\left(\Delta^k\partial_j
g*\frac {1}{|y|^{n-\alpha}}\right)(z)dz\equiv A+B.}
\end{array}
\end{equation}

To deal with $A$ we use that $T*k^j$ is a bounded function. Notice that for $x\in
(3Q)^c$ and $y\in 2Q$ we have

$$|g(y)|\leq \frac{C\|\varphi_Q\|_\infty}{\ell(Q)^\alpha}.$$

Let $Q_0$ stand for the unit cube centered at $0$. Moving $\Delta^k\partial_j$ from
$g$ to $\frac{1}{|y|^{n-\alpha}}$ and making the obvious change of variables one gets

$$\begin{array}{l} \displaystyle{|A|\leq C\sup_{1\leq i\leq
n}\|T*k^i\|_\infty\frac{\|\varphi_Q\|_\infty}{\ell(Q)^\alpha}\int_{(3Q)^c}\int_{2Q}\frac{dydz}{|z-y|^{2n-\alpha}}}\\\\
\displaystyle{\leq C\sup_{1\leq i\leq
n}\|T*k^i\|_\infty\int_{(3Q_0)^c}\int_{2Q_0}\frac{dydz}{|z-y|^{2n-\alpha}}\leq
C\sup_{1\leq i\leq n}\|T*k^i\|_\infty.}
\end{array}$$

Let's now turn our attention to $B$. Recall that we have

\begin{equation}
\label{leibniz} \Delta^k(hg)=\sum_{i_1,...,i_k=1}^n\sum_{l_1,...,l_k=0}^2
\left(\begin{array}{l} 2\\l_1\end{array}\right)...
\left(\begin{array}{l}2\\l_k\end{array}\right)
\partial_{i_1...i_k}^{l_1...l_k}h\;\partial_{i_1...i_k}^{2-l_1...2-l_k}g,
\end{equation}

\noindent where
$\displaystyle{\partial_{i_1...i_k}^{l_1...l_k}=(\partial_{i_1})^{l_1}...(\partial_{i_k})^{l_k}}$.

Since

$$\Delta^k(\partial_j g)=\Delta^k\left(k_x^i\;\partial_j\varphi_Q
\right)+\Delta^k\left(\varphi_Q\;\partial_j k_x^i\right),$$

we have

\begin{equation}
\label{casB}
\begin{array}{l}
\displaystyle{B\leq C\sup_{1\leq i\leq
n}\|T*k^i\|_\infty\int_{3Q}\left|\left(\Delta^k\left(k_x^i\;\partial_j\varphi_Q
\right)*\frac 1 {|y|^{n-\alpha}}\right)(z)\right|dz}\\\\
\displaystyle{+C\sup_{1\leq i\leq
n}\|T*k^i\|_\infty\int_{3Q}\left|\left(\Delta^k\left(\varphi_Q\;\partial_j
k_x^i\right)*\frac 1 {|y|^{n-\alpha}}\right)(z)\right|dz}\\\\
\displaystyle{\equiv C\sup_{1\leq i\leq
n}\|T*\frac{x_i}{|x|^{1+\alpha}}\|_\infty\left(B_1+B_2\right).}
\end{array}
\end{equation}

Using (\ref{leibniz}), support $\varphi_Q \subset 2Q,$
$\|\partial^s\varphi_Q\|_\infty\leq C_s\ell(Q)^{-|s|}$, $|s|\geq 0$, $x\notin 3Q$ and
changing variables, we get

$$\begin{array}{l}
\displaystyle{B_1\leq\sum_{i_1,...,i_k=1}^n
\sum_{l_1,...,l_k=0}^2\frac{C}{\ell(Q)^{l_1+...+l_k+1}}
\int_{3Q}\int_{2Q}\frac{dzdy}{|z-y|^{n-\alpha}|x-y|^{\alpha+2-l_1+...+2-l_k}}}\\\\
\displaystyle{\leq\frac{C}{\ell(Q)^{n+\alpha}}\int_{3Q}\int_{2Q}\frac{dzdy}{|z-y|^{n-\alpha}}=
\frac{C\ell(Q)^{2n}}{\ell(Q)^{n+\alpha+n-\alpha}}
\int_{3Q_0}\int_{2Q_0}\frac{dzdy}{|z-y|^{n-\alpha}}}\\\\
\displaystyle {\leq C.}
\end{array}$$

Arguing similarly we obtain $B_2\leq C$ and therefore we conclude that

$$A+B\leq C\sup_{1\leq i\leq n}\|T*k^i\|_\infty.$$

\item[Case 2:] $x\in 3Q$. Without loss of generality assume $x=0$. Now the function
$g(y)=-\varphi_Q(y)k^i(y)$ may not be smooth, but (\ref{representationformula}) still
holds in the distributions sense. In fact, a different version of
(\ref{representationformula}) will be used for this case. Since $\alpha$ is
non-integer and $[\alpha]=2d$ we readily get
\begin{equation}
\label{representationformula2} f=C\sum_{j=1}^n\Delta^{k-d}\partial_j f*\frac
1{|x|^{n-\alpha+2d}}*\frac{x_j}{|x|^{1+\alpha}},
\end{equation}
where $C=C(n,\alpha)$ and the above identity holds in the distributions sense.

Define $f=T*\frac{1}{|x|^{\alpha-1}}.$ Since $\partial_jf=C(T*k^j)$ and the $T*k^j$
are bounded, the function $f$ satisfies a Lipschitz condition of order $1.$ We get

$$\begin{array}{l}
\displaystyle{\left(\varphi_QT*k^i\right)(0)=<T,g>=c_{n,\alpha}
\sum_{j=1}^n<T*k^j,\Delta^{k-d}\partial_j
g*\frac {1}{|y|^{n-\alpha+2d}}>}\\\\
\displaystyle{=C\sum_{j=1}^n<\partial_j\left(f-f(0)\right)
,\Delta^{k-d}\partial_j g*\frac{1}{|y|^{n-\alpha+2d}}>}.\\\\
\end{array}$$

We claim now that integrating by parts gives

\begin{equation}
\label{ultimpatiment}
\begin{array}{l}
\displaystyle{\sum_{j=1}^n<\partial_j\left(f-f(0)\right),\Delta^{k-d}\partial_j
g*\frac {1}{|y|^{n-\alpha+2d}}>}\\\\ \displaystyle{=<
f-f(0),\Delta^{k-d+1}g*\frac{1}{|y|^{n-\alpha+2d}}>+O\left(\sup_{1\leq i\leq
n}\|T*k^i\|_\infty\right).}
\end{array}
\end{equation}

We postpone the proof of (\ref{ultimpatiment}) and we continue with the argument. If
(\ref{ultimpatiment}) holds, then we can write

$$\begin{array}{l}
\displaystyle{\left|\left(\varphi_QT*k^i\right)(0)\right|\leq C\left|\int_{(3Q)^c}
\left(f(z)-f(0)\right)\left(\Delta^{k+1-d}
g*\frac {1}{|y|^{n-\alpha+2d}}\right)(z)dz\right|}\\\\
\displaystyle{+C\left|\int_{3Q}\left(f(z)-f(0)\right)\left(\Delta^{k+1-d}
g*\frac{1}{|y|^{n-\alpha+2d}}\right)(z)dz\right| +C\sup_{1\leq i\leq
n}\|T*k^i\|_\infty.}
\end{array}$$

Set

$$A=\int_{(3Q)^c}\left(f(z)-f(0)\right)\left(\Delta^{k+1-d} g*\frac {1}{|y|^{n-\alpha+2d}}\right)(z)dz$$

\noindent and

$$B=\int_{3Q}\left(f(z)-f(0)\right)\left(\Delta^{k+1-d} g*\frac {1}{|y|^{n-\alpha+2d}}\right)(z)dz.$$

Using the boundedness of the function $T*k^j=\partial_j f$, Fubini and changing
variables we obtain

$$\begin{array}{l} \displaystyle{|A|\leq C\sup_{1\leq i\leq
n}\|T*k^i\|_\infty\sum_{j=1}^n\int_{(3Q)^c}
|z|\int_{2Q}\frac{|g(y)|}{|z-y|^{2n+1-\alpha}}dydz}\\\\
\displaystyle{\leq C \sup_{1\leq i\leq
n}\|T*k^i\|_\infty\|\varphi_Q\|_\infty\sum_{j=1}^n\int_{(3Q)^c}\int_{2Q}\frac{|z-y|+|y|}{|y|^\alpha|z-y|^{2n+1-\alpha}}dydz}\\\\
\displaystyle{\leq C \sup_{1\leq i\leq n}\|T*k^i\|_\infty\sum_{j=1}^n
\int_{2Q}\frac{1}{|y|^\alpha}\int_{(3Q)^c}\frac{dz}{|z-y|^{2n-\alpha}}dy}\\\\
\displaystyle{+C\sup_{1\leq i\leq
n}\|T*k^i\|_\infty\sum_{j=1}^n\ell(Q)\int_{2Q}\frac{1}{|y|^\alpha}
\int_{(3Q)^c}\frac{dz}{|z-y|^{2n+1-\alpha}}dy}\\\\
\displaystyle{\leq C\sup_{1\leq i\leq n}\|T*k^i\|_\infty.}
\end{array}$$

For the term $B$, write

$$\begin{array}{l} \displaystyle{|B|=\left|\int_{3Q}
\left(f(z)-f(0)\right)\left(\Delta^{k+1-d}g*\frac{1}{|y|^{n-\alpha+2d}}\right)(z)dz\right|}\\\\
\displaystyle{\leq C\left|\int_{3Q}
\sum_{|r|+|s|=n+1-2d}\left(f(z)-f(0)\right)\left(\left(\partial^r\varphi_Q\partial^sk^i\right)*\frac
1{|y|^{n-\alpha+2d}}\right)(z)dz\right|,}
\end{array}$$

\noindent where the last sum is over those multi-indexes $r$ and $s$ that appear in
distributing between $\varphi_Q$ and $k^i$ the $n+1-2d$ derivatives coming from
$\Delta^{k+1-d}$. We will now divide the above sum in two parts, the first one
containing the indexes $|r|\geq 2$ and the second one the remaining indexes. In order
to be able to estimate the integral of this second part, which is the worse, we will
have to subtract a Taylor polynomial of $\varphi_Q$ of order one. Let

$$R(y)=\varphi_Q(y)-\sum_{|m|=0}^1\partial^m\varphi_Q(0)y^m.$$

Then

\begin{equation}
\label{casos}
\begin{array}{l}
\displaystyle{|B|\leq C\sum_{|r|\geq 2}\int_{3Q}\left|f(z)-f(0)\right|
\int_{2Q}\frac{dydz}{\ell(Q)^{|r|}|y|^{\alpha+n+1-2d-|r|}|z-y|^{n-\alpha+2d}}}\\\\
\displaystyle{+C\int_{3Q}\left|f(z)-f(0)\right|\left|\sum_{\tiny{\begin{array}{l}
|r|+|s|=n+1-2d\\|r|\leq 1\end{array}}} \int\frac{\partial^r
R(y)\partial^sk^i(y)}{|z-y|^{n-\alpha+2d}}dy\right|dz}\\\\
\displaystyle{+C\left|\int_{3Q}\left(f(z)-f(0)\right)\sum_{|m|=
0}^1\partial^m\varphi(0)\left(y^m\Delta^{k+1-d}k^i*\frac{1}{|y|^{n-\alpha+2d}}\right)(z)dz\right|}\\\\
\displaystyle{+C\sup_{|m|=1}|\partial^m\varphi_Q(0)|\left|\int_{3Q}\left(f(z)-f(0)\right)\sum_{|s|=n-2d}
\left(\partial^s k^i*\frac 1{|y|^{n-\alpha+2d}}\right)(z)dz\right|}\\\\
\displaystyle{\equiv B_1+B_2+B_3+B_4.}
\end{array}
\end{equation}

Notice that if $|r|\geq 2,$ then we have $\alpha+n+1-2d-|r|\leq\alpha+n-1-2d<n$. Hence
using the boundedness of the functions $T*k^i$, $1\leq i\leq n$, we conclude that
$B_1$ is finite and, by homogeneity, independent of $\ell(Q)$. Thus, $$B_1\leq
C\sup_{1\leq i\leq n}\|T*k^i\|_\infty.$$

We deal now with $B_2$. Write

$$\begin{array}{l}
\displaystyle{B_2=C\int_{3Q}\left|f(z)-f(0)\right|\left|\sum_{\tiny{\begin{array}{l}
|r|+|s|=n+1-2d\\|r|\leq 1\end{array}}}\int_{4Q}\frac{\partial^r
R(y)\partial^sk^i(y)}{|z-y|^{n-\alpha+2d}}dy\right|dz}\\\\
\displaystyle{+C\int_{3Q}\left|f(z)-f(0)\right|\left|\sum_{\tiny{\begin{array}{l}|r|+|s|=n+1-2d\\|r|\leq
1\end{array}}} \int_{(4Q)^c}\frac{\partial^r
R(y)\partial^sk^i(y)}{|z-y|^{n-\alpha+2d}}dy\right|dz=B_{21}+B_{22}.}
\end{array}$$

For the integral over $4Q$, we have to use the Taylor expansion to get integrability.
For the terms with $|r|=1$ we use that
$$|\partial^rR(y)|=|\partial ^r \varphi_Q(y)-\partial ^r
\varphi_Q(0)|\le C\frac{|y|}{\ell(Q)2}$$ and for the term with $|r|=0$ $$ |R(y)|\le
C\frac{|y|^2}{\ell(Q)2}. $$ Therefore
$$\begin{array}{l} \displaystyle{B_{21}\leq C\sup_{1\leq i\leq
n}\|T*k^i\|_\infty\ell(Q)\int_{3Q}\int_{4Q}\frac{|y|}{\ell(Q)2|y|^{\alpha+n-2d}|z-y|^{n-\alpha+2d}}dydz}\\\\
\displaystyle{+C\sup_{1\leq i\leq n}\|T*k^i\|_\infty\ell(Q)
\int_{3Q}\int_{4Q}\frac{|y|^2}{\ell(Q)2|y|^{\alpha+n-2d-1}|z-y|^{n-\alpha+2d}}dydz.}\\\\
\displaystyle{\leq C\sup_{1\leq i\leq
n}\|T*k^i\|_\infty\ell(Q)^{-1}\int_{3Q}\int_{4Q}\frac{dydz}{|y|^{\alpha+n-2d-1}|z-y|^{n-\alpha+2d}}.}
\end{array}$$

Then by homogeneity and local integrability,

$$B_{21}\leq C\sup_{1\leq i\leq n}\|T*k^i\|_\infty.$$

For the integral over $(4Q)^c$, we do not apply Taylor's formula; we just estimate
term by term. For $|r|=0$ (and then $|s|=n+1-2d$) we have that

$$\left|R(y)\partial^sk^i(y)\right|\leq
\frac{C|y|}{\ell(Q)|y|^{\alpha+n+1-2d}}=\frac{C}{\ell(Q)|y|^{n+\alpha-2d}}.$$

For $|r|=1$ the term $\left|\partial^rR(y)\partial^sk^i(y)\right|$ can be estimated by
$C\ell(Q)^{-1}|y|^{-\alpha-n+2d}$, because now $|s|=n-2d$. Therefore

$$\begin{array}{l} \displaystyle{B_{22}\leq C\sup_{1\leq i\leq
n}\|T*k^i\|_\infty\ell(Q)
\int_{3Q}\int_{(4Q)^c}\frac{dy}{\ell(Q)|y|^{\alpha+n-2d}|z-y|^{n-\alpha+2d}}dz}\\\\
\displaystyle{\leq C\sup_{1\leq i\leq n}\|T*k^i\|_\infty.}
\end{array}$$

For $B_3$, separate the terms according to whether $|m|=0$ or $|m|=1$ as follows:

$$\begin{array}{l}
\displaystyle{B_3=\left|\int_{3Q}\left(f(z)-f(0)\right)\varphi_Q(0)\left(\Delta^{k+1-d}k^i*\frac{1}{|y|^{n-\alpha+2d}}\right)(z)dz\right|}\\\\
\displaystyle{+\left|\int_{3Q}\left(f(z)-f(0)\right)\sum_{|m|=1}\partial^m\varphi_Q(0)
\left(y^m\Delta^{k+1-d}k^i*\frac{1}{|y|^{n-\alpha+2d}}\right)(z)dz\right|}\\\\
\displaystyle{\equiv B_{31}+B_{32}.}
\end{array}$$

Now we treat the term $B_{31}$. Taking Fourier transforms on the convolution
$\Delta^{k+1-d}k^i*\frac{1}{|y|^{n-\alpha+2d}}$ we obtain for an appropriate constant
$C$,

$$\widehat{\left(\Delta^{k+1-d}k^i*\frac{1}{|y|^{n-\alpha+2d}}\right)}(\xi)=C\xi_i.$$

Thus

$$\Delta^{k+1-d}k^i*\frac{1}{|y|^{n-\alpha+2d}}=C\partial_i\delta.$$

Hence, by a standard regularization process that we omit,

$$\begin{aligned}
B_{31}&=C|\varphi_Q(0)<\partial_i\delta,f(z)-f(0)>|=C|\varphi_Q(0)<\delta,\partial_i
f>|=C|\varphi_Q(0)\partial_i f(0)|\\ &\leq
 C\|\varphi_Q\|_{\infty}\sup_{1\leq i\leq n}\|T*k^i\|_\infty.
\end{aligned}$$

To estimate $B_{32}$, we take the Fourier transform of
$y^m\Delta^{k+1-d}k^i*\frac{1}{|y|^{n-\alpha+2d}}$, $|m|=1$. We get

$$\widehat{\left(y^m\Delta^{k+1-d}k^i*\frac{1}{|y|^{n-\alpha+2d}}\right)}(\xi)=
C\partial^m\left(\frac{|\xi|^{2k+2-2d}\xi_i}{|\xi|^{1+n-\alpha}}\right)\frac{1}{|\xi|^{\alpha-2d}}=
C\delta_{m,m_i}+C\frac{\xi^m\xi_i}{|\xi|^{2}},$$

\noindent where $m_i$ is the multi-index with all entries equal to $0$ except the
$i-$th entry which is $1$; $\delta_{m,m_i}$ equals one when $m=m_i$ and zero
otherwise.  Hence

$$y^m\Delta^{k+1-d}k^i*\frac{1}{|y|^{n-\alpha+2d}}=C\delta_{m,m_i}\delta+C
 \mbox{ P.V. }\,\frac{z^mz_i}{|z|^{n+2}},$$

where $P.V.$ stands for principal value. Since $|m|=1$,

$$\begin{aligned}
B_{32}&=\left|\sum_{|m|=1}\partial^m\varphi_Q(0)\int_{3Q}(f(z)-f(0))\frac{z^mz_i}{|z|^{n+2}}\right|
\leq\frac{C}{\ell(Q)}\sup_{1\leq i\leq
n}\|T*k^i\|_\infty\int_{3Q}\frac{dz}{|z|^{n-1}}\\ &\leq C\sup_{ 1\leq
i\leq n}\|T*k^i\|_\infty.\\
\end{aligned}$$

Now we are left with term $B_4$. Taking Fourier transforms on the convolution
$\partial^sk^i*\frac{1}{|y|^{n-\alpha+2d}}$, we obtain

$$\widehat{\left(\partial^sk^i*\frac{1}{|y|^{n-\alpha+2d}}\right)}(\xi)
=C\xi^s\frac{\xi_i}{|\xi|^{1+n-\alpha}}\frac{1}{|\xi|^{\alpha-2d}}=C\frac{\xi^s\xi_i}{|\xi|^{n+1-2d}}.$$

Hence, since $|s|=n-2d,$

$$\partial^sk^i*\frac{1}{|y|^{n-\alpha}}=C\, P.V. \,
\frac{z^sz_i}{|z|^{2|s|+1+2d}}.$$

Arguing as before

$$B_4\leq C\sup_{1\leq i\leq
n}\|T*k^i\|_\infty\frac{1}{\ell(Q)}\int_{3Q}\frac{dz}{|z|^{n-1}}\leq C\sup_{1\leq
i\leq n}\|T*k^i\|_\infty.$$

We still have to show claim (\ref{ultimpatiment}). Let $1\leq j\leq n$ and set

$$\omega_j=(-1)^{j-1}dy_1\wedge\cdots\wedge \widehat{dy_j}\wedge\cdots\wedge dy_n.$$

Then, the Green-Stokes Theorem gives
$$\begin{array}{l}
\displaystyle{\sum_{j=1}^n<\partial_j\left(f-f(0)\right),\Delta^{k-d}\partial_j
g*\frac {1}{|y|^{n-\alpha+2d}}>}\\\\
\displaystyle{=-<f-f(0),\Delta^{k-d+1}g*\frac{1}{|y|^{n-\alpha+2d}}>}\\\\
\displaystyle{+\sum_{j=1}^n\lim_{\ep\to
0}\int_{|y|=\ep^{-1}}(f(y)-f(0))\left(g*\Delta^{k-d}\partial_j\frac{1}{|y|^{n-\alpha+2d}}\right)(y)\;\omega_j}\\\\
\displaystyle{-\sum_{j=1}^n\lim_{\ep\to
0}\int_{|y|=\ep}(f(y)-f(0))\left(\Delta^{k-d}\partial_j
g*\frac{1}{|y|^{n-\alpha+2d}}\right)(y)\;\omega_j.}
\end{array}$$
The integral over the sphere of radius $\epsilon^{-1}$ can be easily estimated by a
constant times $\ep^{n-\alpha}$. Thus we are only left with the integral over the
sphere of radius $\epsilon$. For $1\leq j\leq n$ and for a suitable constant $C$ we
can write (recall that for some constant $C$ depending on $n$ and $\alpha$,
$\frac{C}{|y|^{n-\alpha+2d}}=\Delta\left(\frac{1}{|y|^{n-\alpha+2d-2}}\right)$)
$$\begin{array}{l}
\displaystyle{\int_{|y|=\ep}(f(y)-f(0))
\left(\Delta^{k-d}\partial_j g*\frac{1}{|y|^{n-\alpha+2d}}\right)(y)\;\omega_j}\\\\
\displaystyle{=C\sum_{l=0}^n\int_{|y|=\ep}(f(y)-f(0))\left(\Delta^{k-d}\partial_j\partial_l
g*\frac{y_l}{|y|^{n-\alpha+2d}}\right)(y)\;\omega_j.}
\end{array}$$
Notice that when looking at the above integral, the worst case arises when all the
derivatives $\Delta^{k-d}\partial_j\partial_l$ of the product $g=-\varphi_Q k^i$ are
taken on the factor $k^i$. We will only be concerned with this case. For the other
cases argue like in (\ref{casos}). Recall that
$R(y)=\varphi_Q(y)-\sum_{|m|=0}^1\partial^m\varphi_Q(0)y^m$. To get integrability we
use the Taylor expansion of $\varphi_Q$ up to order $1$. Then for $1\leq j\leq n$ we
have
$$\begin{array}{l}
\displaystyle{\int_{|y|=\ep}(f(y)-f(0))\left(\sum_{l}\varphi_Q
\,\Delta^{k-d}\partial_j\partial_l
k^i*\frac{y_l}{|y|^{n-\alpha+2d}}\right)(y)\;\omega_j}\\\\
\displaystyle{=\int_{|y|=\ep}(f(y)-f(0))\left(\sum_lR\;\Delta^{k-d}\partial_j\partial_l
k^i*\frac{y_l}{|y|^{n-\alpha+2d}}\right)(y)\;\omega_j}\\\\
\displaystyle{+\varphi_Q(0)\int_{|y|=\ep}(f(y)-f(0))\left(\sum_l\Delta^{k-d}\partial_j\partial_l
k^i*\frac{y_l}{|y|^{n-\alpha+2d}}\right)(y)\;\omega_j}\\\\
\displaystyle{+\sum_{|m|=1}\partial^m\varphi_Q(0)\int_{|y|=\ep}(f(y)-f(0))\sum_l\left(y^m\Delta^{k-d}\partial_j\partial_l
k^i*\frac{y_l}{|y|^{n-\alpha+2d}}\right)(y)\;\omega_j}\\\\
\displaystyle{=A_1+A_2+A_3.}
\end{array}$$

We will now show that $A_1$ and $A_3$ converge to zero when $\ep\to 0$ and that $A_2$
is bounded above by $C\sup_{1\leq i\leq n}\|T*k^i\|_\infty$.

For $A_1$ we break the convolution integral into two terms corresponding to $3Q$ and
$(3Q)^c:$

$$\begin{array}{l}
\displaystyle{A_1=\int_{|y|=\ep}(f(y)-f(0))\int_{3Q}R(z)\sum_l\Delta^{k-d}\partial_j\partial_l
k^i(z)\frac{y_l-z_l}{|y-z|^{n-\alpha+2d}}dz\;\omega_j}\\\\
\displaystyle{+\int_{|y|=\ep}(f(y)-f(0))\int_{(3Q)^c}R(z)\sum_l\Delta^{k-d}\partial_j\partial_l
k^i(z)\frac{y_l-z_l}{|y-z|^{n-\alpha+2d}}dz\;\omega_j}\\\\
\displaystyle{=A_{11}+A_{12}.}
\end{array}$$

We deal first with $A_{11}$. Since $|R(z)|\leq C|z|^2\ell(Q)^{-2},$ the product
$R\;\Delta^{k-d}\partial_j\partial_l k^i$ is a locally integrable function. Thus,
using the boundedness of $T*k^j,$ we get

$$|A_{11}|\leq C\ep\int_{|y|=\ep} \int_{3Q}
\frac{dz}{|z|^{n-1+\alpha-2d}|z-y|^{n-\alpha+2d-1}}|\;\omega_j| \leq
C\ep^{1+\alpha-2d}.$$

Since we also have $|R(z)|\leq c|z|\ell(Q)^{-1},$ we obtain

$$|A_{12}|\leq C\ep\int_{|y|=\ep}\int_{(3Q)^c}\frac{dz}{|z|^{n+\alpha-2d}|z-y|^{n-\alpha+2d-1}}|\omega_j|\leq C\ep^n.$$

Thus $A_1$ tends to zero with $\ep$.

To estimate $A_2$, take the Fourier transforms of
$\sum_l\Delta^{k-d}\partial_j\partial_lk^i*\frac{y_l}{|y|^{n-\alpha+2d}}$. Then for an
appropriate constant $C$ one has

$$\begin{array}{l}
\displaystyle{\widehat{\left(\sum_l\Delta^{k-d}\partial_j\partial_lk^i*\frac{y_l}{|y|^{n-\alpha+2d}}\right)}(\xi)}
\displaystyle{=C\sum_l|\xi|^{2k-2d}\xi_j\xi_l\frac{\xi_i}{|\xi|^{1+n-\alpha}}\frac{\xi_l}{|\xi|^{2+\alpha-2d}}}=C\frac{\xi_i\xi_j}{|\xi|^2}.\\\\
\end{array}$$

Thus

$$\sum_l\Delta^{k-d}\partial_j\partial_lk^i*\frac{y_l}{|y|^{n-\alpha+2d}}
=C\,\mbox{ P.V. } \frac{y_iy_j}{|y|^{n+2}}.$$

Hence

$$\begin{array}{l}
\displaystyle{|A_2|=\left|C\varphi_Q(0)\int_{|y|=\ep}(f(y)-f(0))\frac{y_iy_j}{|y|^{n+2}}\omega_j\right|}\\\\
\displaystyle{\leq C\sup_{1\leq i\leq
n}\|T*k^i\|_{\infty}\ep^{1-n}\int_{|y|=\ep}|\omega_j|=C\sup_{1\leq i\leq
n}\|T*k^i\|_{\infty}}.
\end{array}$$

For the last term $A_3$, taking the Fourier transform of $\sum_l
y^m\Delta^{k-d}\partial_j\partial_l k^i*\frac{y_l}{|y|^{n-\alpha+2d}}$, we get for a
suitable constants $C_1$ and $C_2$

$$\begin{array}{l} \displaystyle{\widehat{\left(\sum_l
y^m\Delta^{k-d}\partial_j\partial_l
k^i*\frac{y_l}{|y|^{n-\alpha+2d}}\right)}(\xi)}\\\\
\displaystyle{=C_1\sum_l\partial^m\left(\frac{|\xi|^{2k-2d}\xi_j\xi_l\xi_i}{|\xi|^{1+n-\alpha}}\right)\frac{\xi_l}{|\xi|^{2+\alpha-2d}}
=C_1\sum_l\partial^m\left(\frac{\xi_j\xi_l\xi_i}{|\xi|^{2+2d-\alpha}}\right)\frac{\xi_l}{|\xi|^{2+\alpha-2d}}}\\\\
\displaystyle{=C_1\sum_l\left(\delta_{m,m_j}\frac{\xi_i\xi_l}{|\xi|^{2+2d-\alpha}}+\delta_{m,m_i}\frac{\xi_j\xi_l}{|\xi|^{2+2d-\alpha}}+
\delta_{m,m_l}\frac{\xi_i\xi_j}{|\xi|^{2+2d-\alpha}}+C_2\frac{\xi_j\xi_i\xi_l\xi^m}{|\xi|^{4+2d-\alpha}}\right)\frac{\xi_l}{|\xi|^{\alpha+2-2d}}}\\\\
\displaystyle{=C_1\left(\delta_{m,m_j}\frac{\xi_i}{|\xi|^2}+\delta_{m,m_i}\frac{\xi_j}{|\xi|^2}+C_2\frac{\xi^m\xi_j\xi_i}{|\xi|^4}+
\sum_l\delta_{m,m_l}\frac{\xi_i\xi_j\xi_l}{|\xi|^4}\right).}
\end{array}$$

Hence

$$\begin{array}{l} \displaystyle{\left(\sum_l
y^m\Delta^{k-d}\partial_j\partial_l
k^i*\frac{y_l}{|y|^{n-\alpha+2d}}\right)(y)}\\\\
\displaystyle{=C_3\left(\delta_{m,m_j}\frac{y_i}{|y|^n}+\delta_{m,m_i}\frac{y_j}{|y|^n}+C_2\frac{y^my_jy_i}{|y|^{n+2}}+
\sum_l\delta_{m,m_l}\frac{y_iy_jy_l}{|y|^{n+2}}\right)}
\end{array}$$

and since $|m|=1$,

$$\begin{array}{l}
\displaystyle{|A_3|\leq C\int_{|y|=\ep}\frac{|f(y)-f(0)|}{|y|^{n-1}}|\omega_j| \leq
C\ep^{2-n}\int_{|y|=\ep}|\omega_j|=C\ep,}
\end{array}$$

\noindent which completes the proof of claim (\ref{ultimpatiment}).
\end{enumerate}

If $n$ is odd, $\alpha$ non-integer, but $[\alpha]=2d+1$ is also odd, then we replace
the reproducing formula (\ref{representationformula2}) by
\begin{equation}
\label{rp3}
f=C\sum_{j=1}^n\Delta^{k-d}f*\frac{x_j}{|x|^{n-\alpha+2d+2}}*\frac{x_j}{|x|^{1+\alpha}}.
\end{equation}

If $n$ is odd and $\alpha$ is an odd integer of the form $\alpha=2d+1,$ then we use
the reproducing formula (\ref{representationformula2}). Instead of applying Taylor's
expansion up to order $1$, we need in this case to apply Taylor's expansion up to
order $2.$

If $n$ is odd and $\alpha$ is an even integer of the form $\alpha=2d+2$ we use again
formula (\ref{rp3}) and Taylor's expansion up to order $2.$

If $n$ is even we use suitable reproducing formulas (see Lemma 11 in \cite{laura}) and
Taylor's expansions up to order $1$ if $\alpha$ is non-integer and up to order 2 if
$\alpha$ is integer.

\qed\newline

\section{Proof of the Theorem.}

Let $\mu$ be a positive Radon measure and $0<\alpha<1$. For $x\in\Rn$, set

$$\begin{array}{l}
\displaystyle{p_\alpha2(\mu)(x)=\int_{\Rn}\;\int_{\Rn} p_\alpha(x,y,z)d\mu(y)d\mu(z),}\\\\
\displaystyle{M_\alpha\mu(x)=\sup_{r>0}\frac{\mu(B(x,r))}{r^\alpha}}
\end{array}$$

\noindent and

$$U_\alpha^\mu(x)=M_\alpha\mu(x)+p_\alpha(\mu)(x).$$

Recall that in section 2 we defined $p_\alpha(\mu)=\int_{\Rn}\int_{\Rn}\int_{\Rn}
p_\alpha(x,y,z)d\mu(x)d\mu(y)d\mu(z).$ Observe that $p_\alpha(\mu)=\int_{\Rn}
p_\alpha2(\mu)(x)d\mu(x)$. $U_\alpha^\mu$ is the analogue of the potential introduced
in \cite{t1}. The energy associated to this potential is

$$E_\alpha(\mu)=\int_{\Rn} U_\alpha^\mu(x)d\mu(x).\vspace{.5cm}$$

\begin{lemma}
\label{gamma+teo} For each compact set $K\subset\Rn$ and $0<\alpha<1$ we have

$$\gamma_{\alpha,+}(K)\approx\sup_\nu\frac 1 {E_\alpha(\nu)},$$

\noindent where the supremum is taken over the probability measures $\nu$ supported on
$K$.
\end{lemma}

{\em Proof .} Take a positive Radon measure $\mu$ supported on $K$ such that
$\left|\left(\frac{x_i}{|x|^{1+\alpha}}*\mu\right)(x)\right|\leq 1$ for almost all
$x\in\Rn$, $1\leq i\leq n$. We claim that

$$\mu(B(x,r))\leq C r^\alpha,\;\;x\in\Rn,\;\;r>0.$$

To prove the claim take an infinitely differentiable function $\varphi$, supported on
$B(x,2r)$ such that $\varphi=1$ on $B(x,r)$, and $\|\partial^s\varphi\|_\infty\leq C_s
r^{-|s|}$, $|s|\geq 0$. Assume first that $n$ is odd and of the form $n=2k+1$. Then,
by Lemma 11 in \cite{laura},

$$\begin{array}{l} \displaystyle{\mu(B(x,r))\leq\int\varphi
d\mu=c_{n,\alpha}\int\left(\sum_{i=1}^n\Delta^k\partial_i\varphi*\frac
1{|x|^{n-\alpha}}*\frac{x_i}{|x|^{1+\alpha}}\right)(y)d\mu(y)}\\\\
\displaystyle{=-c_{n,\alpha}\sum_{i=1}^n\int\left(\mu*\frac{x_i}{|x|^{1+\alpha}}\right)(y)\left(\Delta^k\partial_i\varphi*\frac
1{|x|^{n-\alpha}} \right)(y)dy}\\\\ \displaystyle{\leq
C\sum_{i=1}^n\left(\int_{B(x,3r)}\left|\left(\Delta^k\partial_i\varphi*\frac
1{|x|^{n-\alpha}}\right)(y)\right|dy+\int_{\Rn\setminus
B(x,3r)}\left|\left(\Delta^k\partial_i\varphi*\frac
1{|x|^{n-\alpha}}\right)(y)\right|dy\right).}
\end{array}$$

Arguing as in Lemma \ref{localitzacio} we get that the last two integrals can be
estimated by $Cr^\alpha$.

If $n$ is even we use the corresponding representation formula in Lemma 11 of
\cite{laura}.

On the other hand, it can be easily shown that

$$\left|R_{\alpha,\ep}(\mu)(x)\right|\leq C,\;\;x\in\Rn,\;\ep>0,$$

\noindent and so, by (\ref{inequa}), we obtain

$$p_\alpha(\mu)\leq C\|\mu\|.$$

By Schwarz inequality

$$E_\alpha(\mu)\leq C\|\mu\|+\|\mu\|^{1/2}p_\alpha(\mu)^{1/2}\leq C\|\mu\|.\vspace{.5cm}$$

Set $\nu=\mu/\|\mu\|$, so that

$$E_\alpha(\nu)=\frac{E_\alpha(\mu)}{\|\mu\|^2}\leq\frac C{\|\mu\|},$$

\noindent and consequently

$$\gamma_{\alpha,+}(K)\leq C\sup_\nu\frac 1{E_\alpha(\nu)}.\vspace{.5cm}$$

The reverse inequality is proved as in \cite{t1} and involves the $T(1)$-Theorem for
non-doubling measures.\qed\newline

It is a crucial fact that the capacity $C_{s,p}$ can be described by means of Wolff
potentials. The Wolff potential of a positive Radon measure $\mu$ is defined by

$$W^\mu(x)=W_{s,p}^\mu(x)=\int_0^\infty\left(\frac{\mu(B(x,r))}{r^{n-sp}}\right)^{p'-1}\frac{dr}r,\;\;x\in\Rn,$$

\noindent where $p'=p/(p-1)$ is the exponent conjugate to $p$.

The Wolff energy of $\mu$ is

$$E(\mu)=E_{s,p}(\mu)=\int_{\Rn} W^\mu(x)d\mu(x).$$

By Wolff's inequality (\cite{adamshedbergllibre}, Theorem 4.5.4, p.110) one has

$$C^{-1}\sup_\mu\frac 1 {E_{s,p}(\mu)^{p-1}}\leq C_{s,p}(K)\leq
C\sup_\mu\frac 1 {E_{s,p}(\mu)^{p-1}},$$

\noindent where $C$ is a positive constant depending only on $s$, $p$ and $n$, and the
supremum is taken over the probability measures $\mu$ supported on $K$.
\newline

\begin{lemma}
\label{equivalencia} For each positive finite Radon measure $\mu$ on $\Rn$ we have

$$p_{\alpha}(\mu)\approx E_{\frac 2 3(n-\alpha),\frac 3 2}(\mu)=\int_{\Rn}\;\int_0^{\infty}
\left(\frac{\mu(B(x,r))}{r^{\alpha}}\right)2\frac{dr}{r}d\mu(x).$$
\end{lemma}

{\em Proof.} Suppose that

$$\int_{\Rn}\;\int_0^{\infty}\left(\frac{\mu(B(x,r))}{r^{\alpha}}\right)2\frac{dr}{r}d\mu(x)<\infty$$

\noindent and set $G=\{(x_1,x_2,x_3):|x_1-x_2|\leq |x_1-x_3|\leq |x_2-x_3|\}$. Using
Lemma \ref{acotaci} and Riemann-Stieltjes integration, we obtain

\begin{equation}
\label{ultima}
\begin{array}{l}
\displaystyle{p_{\alpha}(\mu)=3\underset{G}{\iiint} p_{\alpha}(x_1,x_2,x_3)
d\mu(x_1)d\mu(x_2)d\mu(x_3)}\\\\
\displaystyle{\approx\iint\underset{B(x_3,|x_2-x_3|)}{\int} |x_2-x_3|^{-2\alpha}d\mu(x_1)d\mu(x_2)d\mu(x_3)}\\\\
\displaystyle{=\int_{\Rn}\!\int_{\Rn}\!\frac{\mu(B(x_3,|x_2-x_3|))}{ |x_2-x_3|^{2\alpha}} d\mu(x_2)d\mu(x_3)}\\\\
\displaystyle{=\int_{\Rn}\!\int_0^{\infty}\frac{\mu(B(x_3,r))}{r^{2\alpha}}d\mu(B(x_3,r))d\mu(x_3).}
\end{array}
\end{equation}

Notice that

\begin{equation}
\label{zerolim1} \lim_{r\to\infty}\left(\frac{\mu(B(x,r))}{r^{\alpha}}\right)2
\leq\lim_{r\to\infty}\left(\frac{\mu(\Rn)}{r^{\alpha}}\right)2=0.
\end{equation}

Moreover,

$$\int_\rho^{2\rho}\left(\frac{\mu(B(x,r))}{r^{\alpha}}\right)2\frac{dr}{r}\geq \mu(B(x,\rho))2
\int_\rho^{2\rho}\frac{dr}{r^{2\alpha+1}}=C\left(\frac{\mu(B(x,\rho))}{\rho^{\alpha}}\right)2.\vspace{.5cm}$$

Thus

\begin{equation}
\label{zerolim2} \lim_{r\rightarrow
0}\left(\frac{\mu(B(x,r))}{r^{\alpha}}\right)2=0.\vspace{.5cm}
\end{equation}

Integration by parts in the last integral of (\ref{ultima}), together with
(\ref{zerolim1}) and (\ref{zerolim2}), show that

$$p_\alpha(\mu)\approx\int_{\Rn}\int_0^\infty\left(\frac{\mu(B(x,r))}{r^\alpha}\right)2\frac{dr}rd\mu(x).$$

Suppose now that $p_\alpha(\mu)<\infty$. We claim that we can assume that

\begin{equation}
\label{liminf} \underset{r\to 0}{\lim\inf}\;\frac{\mu(B(x,r))}{r^\alpha}=0,\;\mbox{
for $\mu$-almost all }\;x\in\Rn.
\end{equation}

If (\ref{liminf}) holds, then integrating by parts in the last integral of
(\ref{ultima}) one can deduce that

$$p_\alpha(\mu)\approx\int_{\Rn}\;\int_0^{\infty}\left(\frac{\mu(B(x,r))}{r^{\alpha}}\right)2\frac{dr}{r}d\mu(x),$$

\noindent and in this case we are done.\newline

Otherwise there exists a $\mu$-measurable set $F$ such that $\mu(F)>0$ and

$$\underset{r\to 0}{\lim\inf}\;\frac{\mu(B(x,r))}{r^\alpha}>0,\;\;x\in F.$$

Shrinking $F$ we can assume that

$$\underset{r\to 0}{\lim\inf}\;\frac{\mu(B(x,r))}{r^\alpha}>a>0,\;\;x\in F.$$

By Egorov we can find $r_0>0$ and a $\mu$-measurable subset $G$ of $F$ such that
$\mu(G)>0$ and

\begin{equation}
\label{boundinf} \mu(B(x,r))>\frac a 2 \;r^{\alpha},\; x\in G\;\mbox{ and } r\leq r_0.
\end{equation}

>From (\ref{ultima}) we get, applying (\ref{boundinf}) twice,

$$\begin{array}{l}
\displaystyle{p_{\alpha}(\mu)\approx\int_{\Rn}\!\int_{\Rn}\!\frac{\mu(B(x_3,|x_2-x_3|))}{
|x_2-x_3|^{2\alpha}} d\mu(x_2)d\mu(x_3)}\\\\
\displaystyle{\geq\int_{G}\;\int_{B(x_3,r_0)}\frac{\mu(B(x_3,|x_2-x_3|))}{|x_2-x_3|^{2\alpha}}d\mu(x_2)d\mu(x_3)}\\\\
\displaystyle{\geq\frac a
2\int_G\;\int_{B(x_3,r_0)}\frac{d\mu(x_2)d\mu(x_3)}{|x_2-x_3|^{\alpha}}}\\\\
\displaystyle{=\frac a 2 \int_G\!\int_0^{\infty}\mu(\{x_2\in
B(x_3,r_0):\;|x_2-x_3|^{-\alpha}\geq t\})dt d\mu(x_3)}\\\\
\displaystyle{\geq\frac{a\alpha}2\int_G\!\int_0^{r_0}\frac{\mu(B(x_3,r))}{r^{1+\alpha}}drd\mu(x_3)}\\\\
\displaystyle{\geq\frac{a2\alpha}2\int_G\!\int_0^{r_0}\frac{dr}r=+\infty,}
\end{array}$$

\noindent which is a contradiction.

\qed\newline

{\bf Remark. } In Theorem 2.2 of \cite{mattilacantor} it is shown that for any finite
Borel measure in $\C$, one has the following inequality,

\begin{equation}
\label{desimattila} \int_\C\;\int_\C\;\int_\C
c2(x_1,x_2,x_3)d\mu(x_1)d\mu(x_2)d\mu(x_3)\leq
C\int_\C\;\int_0^\infty\frac{\mu(B(x,r))2}{r2}\frac{dr}{r}d\mu(x).\vspace{.5cm}
\end{equation}

On the other hand, for $\alpha=1,$ there is no general lower inequality like the one
in Lemma \ref{acotaci}. Although we have

$$c(x_1,x_2,x_3)\leq \frac 2{|x_2-x_3|},$$

\noindent the reverse inequality may fail very badly. Thus the reverse inequality in
(\ref{desimattila}) does not hold for general measures $\mu$. However, see Theorem 2.3
in \cite{mattilacantor} where a related result is shown when $\mu$ is the Hausdorff
measure associated to some measure function $h$, restricted to some Cantor
sets.\newline

We turn now to the proof of the main Theorem.\newline

{\em Proof of the Theorem.} We deal first with the inequality

\begin{equation}
\label{compupper} C_{\frac 2 3 (n-\alpha),\frac 3 2}(K)\leq C\gamma_{\alpha_+}(K).
\end{equation}

Assume that for a probability measure $\mu$ supported on $K$ we have

$$E_{\frac 2 3 (n-\alpha),\frac 3 2}(\mu)=\int_{\Rn}\;\int_0^\infty\left(\frac{\mu(B(x,r))}{r^\alpha}\right)2\frac{dr} r d\mu(x)\equiv E<\infty.$$

Then by Chebyshev, for each $t>0$,

$$\mu\{x\in K:\;\int_0^\infty\left(\frac{\mu(B(x,r))}{r^\alpha}\right)2\frac{dr} r>t\}\leq\frac E t.$$

Taking $t=2E$, we obtain a compact set $F\subset K$ such that

$$\int_0^\infty\left(\frac{\mu(B(x,r))}{r^\alpha}\right)2\frac{dr} r \leq 2E,\;\;x\in F,$$

\noindent and

$$\mu(F)\geq\frac 1 3.$$

If we set $\nu=\mu_{|F}/\mu(F)$, then for some positive constant $C$ depending on
$\alpha$,

\begin{equation}
\label{superfluous}
C\left(\frac{\nu(B(x,\rho))}{\rho^\alpha}\right)2\leq\int_\rho^{2\rho}\left(\frac{\nu(B(x,r))}{r^\alpha}\right)2\frac{dr}
r \leq 18 E,\;\;x\in F.
\end{equation}

To see that $\nu$ satisfies the $\alpha$-growth condition, notice that if $x\notin F$
and $B(x,r)\cap F=\emptyset$, then $\nu(B(x,r))=0$, and if there is some $\xi\in F\cap
B(x,r)$, then due to (\ref{superfluous})

$$\nu(B(x,r))\leq \nu(B(\xi,2r))\leq Cr^\alpha\sqrt E.$$

Hence we have

$$M_\alpha\nu(x)\leq C\sqrt E,\;x\in\Rn.$$

Then by Lemma \ref{equivalencia} and Schwarz inequality we get

$$E_{\alpha}(\nu)=\int_{\Rn} U_{\alpha}^\nu(x)d\nu(x)\leq C\sqrt E+ p_\alpha(\nu)^{1/2}\leq C\sqrt E.$$

\noindent Thus, by Lemma \ref{gamma+teo}, we obtain

$$E^{-1/2}\leq C E_\alpha(\nu)^{-1}\leq C\gamma_{\alpha,+}(K),$$

\noindent which implies (\ref{compupper}).

To see the reverse inequality, let $\mu$ be a probability measure supported on $K$
such that

$$E_\alpha(\mu)=\int_{\Rn} U_{\alpha}^\mu(x)d\mu(x)<\infty.$$

Since

$$E_\alpha(\mu)\geq\int p_\alpha(\mu)(x)d\mu(x),$$

as before, by Chebyshev,

$$\mu\{x\in K:\;p_\alpha(\mu)(x)>t\}\leq\frac{E_\alpha(\mu)}t,\;t>0.$$

Taking $t=2E_\alpha(\mu)$ we find a compact set $F\subset K$ such that

$$p_\alpha(\mu)(x)\leq 2E_\alpha(\mu),\;\mbox{ for }
x\in F, $$

\noindent and

$$\mu(F)\geq\frac 1 3.$$

Set $\nu=\mu_{|F}/\mu(F)$. Then

$$p_\alpha(\nu)=\int_{F} p_\alpha2(\nu)(x)d\nu(x)\leq 36 E_\alpha(\mu)2,$$

and so, by Lemma \ref{equivalencia}

$$E_\alpha(\mu)^{-1}\leq 6 p_\alpha(\nu)^{-1/2}\approx E_{\frac 2 3(n-\alpha),\frac 3 2}(\nu)^{-1/2}
\leq C_{\frac 2 3(n-\alpha),\frac 3 2}(K),$$

\noindent which ends the proof of the Theorem.

\qed

\section*{Acknowledgements}

The authors were partially supported by grants DGI \, MTM2004-00519 and 2001/SGR/00431
(Generalitat de Catalunya).

Joan Mateu and Joan Verdera

Departament de Matem\`{a}tiques

Universitat Aut\`onoma de Barcelona

08193 Bellaterra (Barcelona)

Spain

{\em E-mail address:} mateu@mat.uab.es, jvm@mat.uab.es.\newline

Laura Prat

Departament de Matem\`atica Aplicada i An\`alisi,

Universitat de Barcelona

Gran Via de les Corts Catalanes, 585

08007 Barcelona

Spain\\

{\em E-mail address:} laura@mat.ub.es

\end{document}